\magnification=\magstep1
\input amstex
\UseAMSsymbols
\input pictex
\vsize=24truecm
\hsize=15truecm
\NoBlackBoxes
\parindent=18pt

   \font\rmk=cmr8      

\font\gr=cmbx10  
\font\gross=cmbx10 scaled\magstep1 

\def\pd{\operatorname{pd}}

\def\findim{\operatorname{fin.dim.}}
\def\Findim{\operatorname{Fin.dim.}}
  \def\ss{\ssize }
\def\arr#1#2{\arrow <1.5mm> [0.25,0.75] from #1 to #2}

\def\s{\hfill \square}

\centerline{\gross Are special biserial algebras homologically tame?}
	\bigskip
\centerline{\gr Claus Michael Ringel}
	\bigskip\bigskip

\noindent {\narrower Abstract: \rmk Birge Huisgen-Zimmermann calls a finite dimensional algebra
homologically tame provided the little and the big finitistic
dimension are equal and finite.
The question formulated in the title has been discussed by her
in the paper ``Representation-tame algebras need not be homologically
tame'', by looking for any $\ss r\ge 1$
at a sequence of algebras $\ss\Lambda_m$ with big finitistic dimension $\ss r+m$.
As we will show, also the little finitistic dimension of $\ss\Lambda_m$ is $\ss r+m$. It follows
that contrary to her assertion, all the algebras $\ss\Lambda_m$ are homologically tame.
	\medskip
\noindent
Key words. Finitistic dimension. Special biserial algebras.

\noindent
2010 Math Subject classification. Primary
16E10. Secondary
16G20, 16G60.\par}
	\bigskip

Birge Huisgen-Zimmermann calls a finite dimensional algebra
{\it homologically tame} provided the little and the big finitistic
dimension are equal and finite.
The question formulated in the title has been discussed by her
in the paper [HZ], by looking for any $r\ge 1$
at a sequence of algebras $\Lambda_m$ with big finitistic dimension $r+m$.
She presented a quite surprising infinite-dimensional $\Lambda_m$-module
with projective dimension $r+m$, stressing that related finite-dimensional
modules have infinite projective dimension. Nonetheless, as we will show, there
do exist finite-dimensional $\Lambda_m$-modules with projective dimension $r+m$.
Thus, also the little finitistic dimension of $\Lambda_m$ is $r+m$. It follows
that contrary to her assertion, all the algebras $\Lambda_m$ are homologically tame.
	\medskip
{\bf Notation.} Let $k$ be a field and $r\ge 1$ a fixed natural
number.
We will deal with a sequence of finite-dimensional $k$-algebras
$\Lambda_m$ with $m\ge 0$, with
$\Lambda_{m}$ being a factor algebra of $\Lambda_{m+1}$ for all $m$
(thus, $\Lambda_{m}$-modules can be considered as $\Lambda_{m+1}$-modules)
such that the projective $\Lambda_{m}$-modules are also projective as
$\Lambda_{m+1}$-modules.
The modules to be considered are (not necessarily finitely generated, left)
$\Lambda_m$-modules for some $m$.
Given any module $M$, we denote by $PM$ a projective cover, by $\Omega M$
the first syzygy module and by $\pd M$ the projective dimension of $M$.
If $x$ is a vertex of a quiver $Q$, the corresponding simple representation will also
be denoted by $x$.
	\medskip
{\bf Outline.}
In section 1, we recall the definition of the special biserial algebras $\Lambda_m$
considered in [HZ].
In section 2, we exhibit for any $m$
a finite-dimensional $\Lambda_m$-module $Z_m$ of
projective dimension $r+m.$ Thus the little finitistic dimension
$\findim\Lambda_m$ of $\Lambda_m$ is at least $r+m.$
Section 3 presents a proof of the assertion in [HZ]
that the big finitistic dimension $\Findim\Lambda_m$
of $\Lambda_m$ is at most $r+m.$
Combining these results, we get
$$
 r+m \le \findim\Lambda_m \le \Findim \Lambda_m \le r+m.
$$
It follows that both the little and the big finitistic
dimension are equal to $r+m.$
	\bigskip
{\bf 1. The algebras $\Lambda_m$.}
As we mentioned, we consider the algebras $\Lambda_m$ exhibited in [HZ]. This means,
we deal with the following quiver with relations.
The labels of the vertices are those used in [HZ], but
we denote all arrows just by $\alpha$ or $\beta$,
so that $\alpha\beta = 0 = \beta\alpha$ (and so that at any vertex, at most one $\alpha$-arrow and at most one $\beta$-arrow start, and similarly,
at most one $\alpha$-arrow and at most one $\beta$-arrow end; finally, the arrow $a_1\to d_0$
is an $\alpha$-arrow).
The $\alpha$-arrows are drawn as solid arrows, the $\beta$-arrows are the dashed ones.
$$
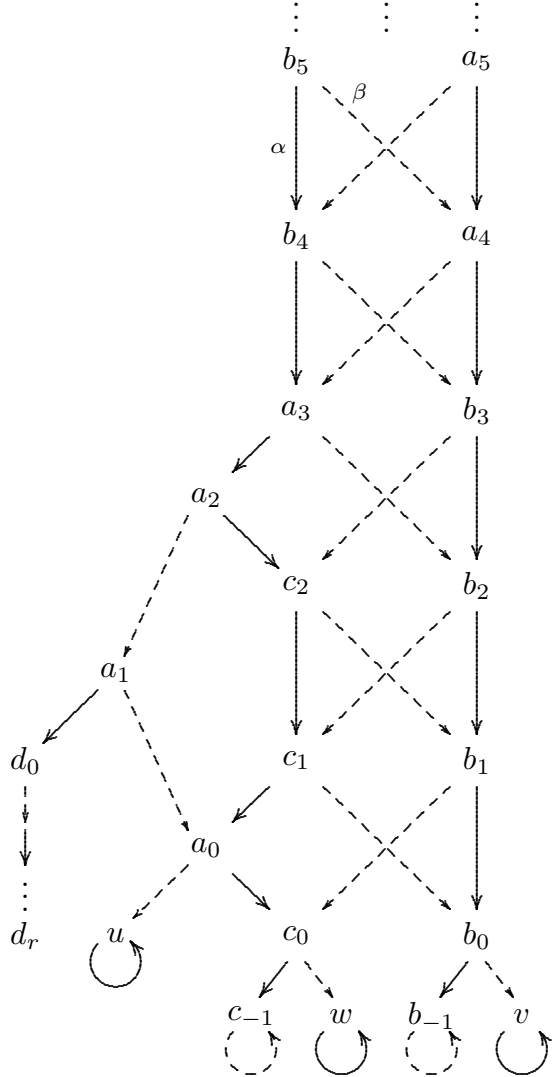
\hbox{\beginpicture
    \setcoordinatesystem units <1cm,.97cm>
\multiput{} at 0 -.5  / 
\put{$d_0$} at 0 3
\put{$d_r$} at 0 1
\put{$u$} at 1 1
\put{$c_{-1}$} at 2.5 0.07
\put{$v$} at 5.5 0.07
\put{$w$} at 3.5 0.07
\put{$b_{-1}$} at 4.5 0.07
\put{$c_0$} at 3 1
\put{$b_0$} at 5 1
\put{$a_0$} at 2 2
\put{$c_1$} at 3 3
\put{$b_1$} at 5 3
\put{$a_1$} at 1 4
\put{$c_2$} at 3 5
\put{$b_2$} at 5 5
\put{$a_2$} at 2 6
\arr{3 4.7}{3 3.3}
\arr{5 4.7}{5 3.3}
\arr{5 2.7}{5 1.3}

\arr{0.8 3.8}{0.2 3,2}
\arr{4.9 0.7}{4.6 0.3}

\arr{2.9 0.7}{2.6 0.3}

\put{$a_3$} at 3 7
\put{$b_3$} at 5 7
\arr{2.7 6.7}{2.3 6.3}
\arr{5 6.7}{5 5.3}

\put{$b_4$} at 3 9
\put{$a_4$} at 5 9
\arr{3 8.7}{3 7.3}
\arr{5 8.7}{5 7.3}

\put{$b_5$} at 3 11
\put{$a_5$} at 5 11
\arr{3 10.7}{3 9.3}
\arr{5 10.7}{5 9.3}

\multiput{$\vdots$} at 3 11.6  4 11.6  5 11.6 /

\multiput{$\ss \alpha$} at  2.8 10   /
\multiput{$\ss \beta$} at  3.7 10.62   /

\circulararc 270 degrees from 3.3 -.1 center at 3.5 -.3
\circulararc 260 degrees from 5.3 -.1 center at 5.5 -.3
\circulararc 280 degrees from .8 .9 center at 1 .7
\arr{1.2 .89}{1.17 .93}
\arr{2.735 -.16}{2.7 -.1}
\arr{3.735 -.16}{3.7 -.1}
\arr{4.735 -.16}{4.7 -.1}
\arr{5.735 -.16}{5.7 -.1}

\arr{0 2.2}{0 1.8}
\put{$\vdots$} at 0 1.55

\arr{2.2 5.8}{2.8 5.2}
\arr{2.3 1.7}{2.7 1.3}
\arr{2.7 2.7}{2.3 2.3}

\setdashes <1mm>
\arr{3.3 6.7}{4.7 5.3}
\arr{4.7 6.7}{3.3 5.3}
\arr{3.3 8.7}{4.7 7.3}
\arr{4.7 8.7}{3.3 7.3}
\arr{3.3 10.7}{4.7 9.3}
\arr{4.7 10.7}{3.3 9.3}

\arr{5.1 0.7}{5.4 0.3}

\arr{3.1 0.7}{3.4 0.3}
\arr{1.8 1.8}{1.2 1,2}
\arr{0 2.7}{0 2.3}
\circulararc 260 degrees from 2.3 -.1 center at 2.5 -.3
\circulararc 270 degrees from 4.3 -.1 center at 4.5 -.3
\arr{3.3 2.7}{4.7 1.3}
\arr{4.7 2.7}{3.3 1.3}
\arr{3.3 4.7}{4.7 3.3}
\arr{4.7 4.7}{3.3 3.3}
\arr{1.8 5.8}{1.1 4.2}
\arr{1.1 3.8}{1.8 2.2}

\endpicture}
$$
The path starting at $d_0$ is an alternating
$\beta$-$\alpha$-path of length $r$ with vertices $d_0, d_1,\dots,d_r.$
The $\alpha$-path of length 2
ending at $a_m$ with $m\ge 3$ starts in $b_{m+2}$. Similarly, the
$\alpha$-path of length 2
ending at $b_m$ with $m\ge 2$ starts in $a_{m+2}$.

There are the following additional relations:
the square of any loop is zero, and we have $\alpha^n = \beta^m$,
whenever this makes sense.

It is easily seen that $\pd d_i = r-i$ (in particular
$\pd d_0 = r$) and that $\pd u = \pd v = \pd w = \pd c_{-1} = \pd b_{-1} =
\infty$.

The algebra $\Lambda_m$ with $m\ge 0$ is given by the full subquiver
with vertices $a_i,b_i,c_i$ where $i \le m$ and all the vertices
$u,v,w,d_0,\dots,d_r$.
		\bigskip
{\bf 2. A finite-dimensional
$\Lambda_m$-module $Z_m$ of projective dimension $r+m$.}
We are going to exhibit a sequence of $\Lambda_m$-modules $Z_m$. All are direct sums
of string modules.
	\vfill\eject
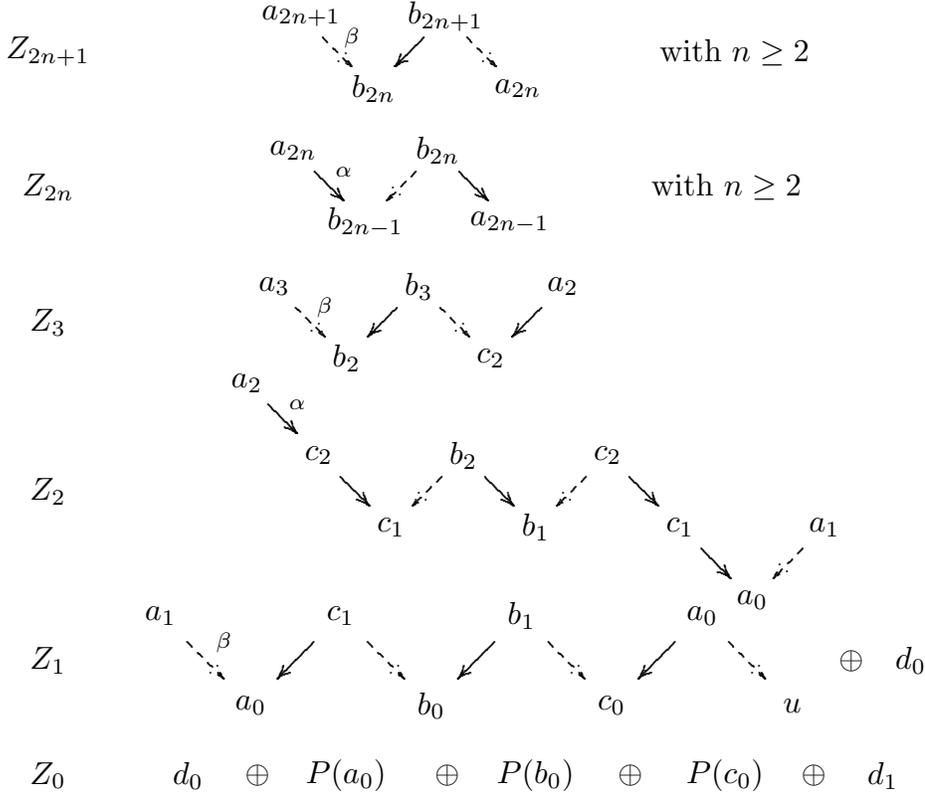
$$
\hbox{\beginpicture
    \setcoordinatesystem units <.9cm,.75cm>
\put{$Z_{2n+1}$} at -6 -0
\put{$Z_{2n}$} at -6 -2
\put{$Z_3$} at -6 -4
\put{$Z_2$} at -6 -6.5
\put{$Z_1$} at -6 -9
\put{$Z_0$} at -6 -10.7
\put{\beginpicture
    \setcoordinatesystem units <.8cm,.8cm>
\put{with $n\ge 2$} at 6 0.5
\put{$a_{2n+1}$} at  0 1
\put{$b_{2n}$} at  1 0
\put{$b_{2n+1}$} at  2 1
\put{$a_{2n}$} at  3 0
\arr{1.7 0.7}{1.3 0.3}
\put{$\ss \beta$} at 0.7 0.7
\setdashes <.8mm>
\arr{0.3 0.7}{0.7 0.3}
\arr{2.3 0.7}{2.7 0.3}

\endpicture} at 0 -0
\put{\beginpicture
    \setcoordinatesystem units <.8cm,.8cm>
\put{with $n\ge 2$} at 6 0.5
\put{$a_{2n}$} at  0 1
\put{$b_{2n-1}$} at  1 0
\put{$b_{2n}$} at  2 1
\put{$a_{2n-1}$} at  3 0
\arr{0.3 0.7}{0.7 0.3}
\arr{2.3 0.7}{2.7 0.3}
\put{$\ss \alpha$} at 0.7 0.7
\setdashes <.8mm>
\arr{1.7 0.7}{1.3 0.3}
\endpicture} at 0 -2
\put{\beginpicture
    \setcoordinatesystem units <.8cm,.8cm>
\put{$a_3$} at  0 1
\put{$b_2$} at  1 0
\put{$b_3$} at  2 1
\put{$c_2$} at  3 0
\put{$a_2$} at  4 1
\arr{1.7 0.7}{1.3 0.3}
\arr{3.7 0.7}{3.3 0.3}
\put{$\ss \beta$} at 0.7 0.7
\put{} at 7 0
\setdashes <.8mm>
\arr{0.3 0.7}{0.7 0.3}
\arr{2.3 0.7}{2.7 0.3}
\endpicture} at -.2 -4
\put{\beginpicture
    \setcoordinatesystem units <.8cm,.8cm>
\put{$a_2$} at  0 3
\put{$c_2$} at  1 2
\put{$c_1$} at  2 1
\put{$b_2$} at  3 2
\put{$b_1$} at  4 1
\put{$c_2$} at  5 2
\put{$c_1$} at  6 1
\put{$a_0$} at  7 0
\put{$a_1$} at  8 1
\arr{0.3 2.7}{0.7 2.3}
\arr{1.3 1.7}{1.7 1.3}
\arr{3.3 1.7}{3.7 1.3}
\arr{5.3 1.7}{5.7 1.3}
\arr{6.3 .7}{6.7 .3}
\put{$\ss \alpha$} at .7 2.7
\setdashes <.8mm>
\arr{2.7 1.7}{2.3 1.3}
\arr{4.7 1.7}{4.3 1.3}
\arr{7.7 .7}{7.3 .3}
\endpicture} at 0 -6.5
\put{\beginpicture
    \setcoordinatesystem units <1cm,1cm>
\put{$a_1$} at  1 2
\put{$a_0$} at  2 1
\put{$c_1$} at  3 2
\put{$b_0$} at  4 1
\put{$b_1$} at  5 2
\put{$c_0$} at  6 1
\put{$a_0$} at  7 2
\put{$u$} at  8 1
\arr{2.7 1.7}{2.3 1.3}
\arr{4.7 1.7}{4.3 1.3}
\arr{6.7 1.7}{6.3 1.3}
\put{$\ss \beta$} at 1.7 1.7
\put{$\oplus \quad d_0$} at 9 1.5
\setdashes <.8mm>
\arr{1.3 1.7}{1.7 1.3}
\arr{3.3 1.7}{3.7 1.3}
\arr{5.3 1.7}{5.7 1.3}
\arr{7.3 1.7}{7.7 1.3}
\endpicture} at 0 -9
\put{\beginpicture
    \setcoordinatesystem units <.7cm,1cm>
\put{$d_0$} at -.5 0
\put{$P(a_0)$} at 2 0
\put{$P(b_0)$} at 5 0
\put{$P(c_0)$} at 8 0
\put{$d_1$} at 10.5 0
\multiput{$\oplus$} at .6 0  3.6 0  6.5 0  9.4 0 /
\endpicture} at 0 -10.7

\endpicture}
$$
	\medskip\smallskip
\noindent
For $Z_m$ with $m$ even, the southeast arrows are $\alpha$-arrows;
for $m$ odd, the southeast arrows are $\beta$-arrows.
	\medskip
{\bf Proposition.} {\it For $m\ge 0$, we have $\Omega Z_{m+1} = Z_{m},$
and $\pd Z_m = r+m$.}
	\medskip
Proof. The first assertion is easily verified.
Since $\pd d_0 = r$ and $\pd d_1 = r-1$, the second assertion is an immediate
consequence, using induction. $\s $
	\smallskip
{\bf Remark.} The modules $Z_m$ with $m\ge 1$ are finite dimensional
$\Lambda_m$-modules, but not $\Lambda_{m-1}$-modules. Since the projective dimension of any
$Z_m$
is finite, the modules $Z_m$ with $m\ge 2$ are counter-examples to Claim 2 of [HZ].
	\bigskip
{\bf 3. The big finitistic dimension of $\Lambda_m.$}
Let $\Lambda_1'$ be obtained from $\Lambda_2$ by deleting the vertices $a_2$ and $b_2.$ We note the following:  {\it Let $M$ be a $\Lambda_m$-module. If
$m =1$ or $m\ge 3$, then $\Omega M$ is a $\Lambda_{m-1}$-module. If $m=2$,
then $\Omega M$ is a $\Lambda_1'$-module.}
	\smallskip
Let $\Cal X$ be the set of the following 10 isomorphism classes of 
$\Lambda_1'$-modules; these are string modules $X$ with $X_{c_2}\neq 0.$
$$
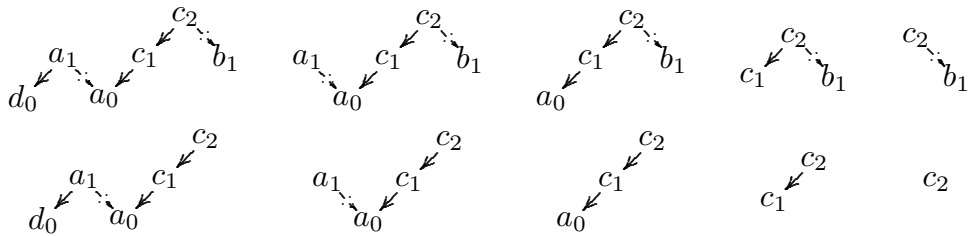
\hbox{\beginpicture
    \setcoordinatesystem units <.45cm,.45cm>
\put{\beginpicture
\put{$d_0$} at 0 0
\put{$a_1$} at 1 1
\put{$a_0$} at 2 0
\put{$c_1$} at 3 1
\put{$c_2$} at 4 2
\put{$b_1$} at 5 1
\arr{0.7 0.7}{0.3 0.3}
\arr{2.7 0.7}{2.3 0.3}
\arr{3.7 1.7}{3.3 1.3}
\setdashes <.8mm>
\arr{1.3 0.7}{1.7 0.3}
\arr{4.3 1.7}{4.7 1.3}
\endpicture} at 0 0
\put{\beginpicture
\put{$a_1$} at 1 1
\put{$a_0$} at 2 0
\put{$c_1$} at 3 1
\put{$c_2$} at 4 2
\put{$b_1$} at 5 1
\arr{2.7 0.7}{2.3 0.3}
\arr{3.7 1.7}{3.3 1.3}
\setdashes <.8mm>
\arr{1.3 0.7}{1.7 0.3}
\arr{4.3 1.7}{4.7 1.3}
\endpicture} at 6.5 0
\put{\beginpicture
\put{$a_0$} at 2 0
\put{$c_1$} at 3 1
\put{$c_2$} at 4 2
\put{$b_1$} at 5 1
\arr{2.7 0.7}{2.3 0.3}
\arr{3.7 1.7}{3.3 1.3}
\setdashes <.8mm>
\arr{4.3 1.7}{4.7 1.3}
\endpicture} at 12 0
\put{\beginpicture
\put{$c_1$} at 3 1
\put{$c_2$} at 4 2
\put{$b_1$} at 5 1
\arr{3.7 1.7}{3.3 1.3}
\setdashes <.8mm>
\arr{4.3 1.7}{4.7 1.3}
\endpicture} at 16.5 0
\put{\beginpicture
\put{$c_2$} at 4 2
\put{$b_1$} at 5 1
\setdashes <.8mm>
\arr{4.3 1.7}{4.7 1.3}
\endpicture} at 20 0
\put{\beginpicture
\put{$d_0$} at 0 0
\put{$a_1$} at 1 1
\put{$a_0$} at 2 0
\put{$c_1$} at 3 1
\put{$c_2$} at 4 2
\arr{0.7 0.7}{0.3 0.3}
\arr{2.7 0.7}{2.3 0.3}
\arr{3.7 1.7}{3.3 1.3}
\setdashes <.8mm>
\arr{1.3 0.7}{1.7 0.3}
\endpicture} at 0 -3
\put{\beginpicture
\put{$a_1$} at 1 1
\put{$a_0$} at 2 0
\put{$c_1$} at 3 1
\put{$c_2$} at 4 2
\arr{2.7 0.7}{2.3 0.3}
\arr{3.7 1.7}{3.3 1.3}
\setdashes <.8mm>
\arr{1.3 0.7}{1.7 0.3}
\endpicture} at 6.5 -3
\put{\beginpicture
\put{$a_0$} at 2 0
\put{$c_1$} at 3 1
\put{$c_2$} at 4 2
\arr{2.7 0.7}{2.3 0.3}
\arr{3.7 1.7}{3.3 1.3}
\endpicture} at 12 -3
\put{\beginpicture
\put{$c_1$} at 3 1
\put{$c_2$} at 4 2
\arr{3.7 1.7}{3.3 1.3}
\endpicture} at 16.5 -3
\put{\beginpicture
\put{$c_2$} at 4 2
\endpicture} at 20 -3

\endpicture}
$$
	\medskip
{\bf Lemma 1.} {\it The modules in $\Cal X$ have infinite projective dimension.}
	\medskip
Proof. For the modules $X$ in the first row,
$c_{-1}$ is a direct summand of $\Omega^2X.$ For the modules $X$ in the second
row, $v$ is a direct summand of $\Omega^3X.$
$\s$
	\medskip
{\bf Lemma 2.} {\it Any $\Lambda_1'$-module
is the direct sum of a $\Lambda_1$-module, of copies of $P(c_2)$, and
of copies of modules in $\Cal X.$}
	\medskip
Proof. Let $M$ be a $\Lambda_1'$-module without a direct summand of the form
$P(c_2)$.
Let $U$ be the subquiver of the quiver of
$\Lambda_1'$ with vertices $d_0, a_1,a_0,c_1,c_2,b_1.$

Since $U$ is a Dynkin quiver,
any representation of $U$ is a direct sum of finite-dimensional indecomposable
representations. We decompose the restriction
$M|U$ of $M$ to $U$ as follows: $M|U = X \oplus Y$, where $X$ is a direct sum of
copies of modules in $\Cal X$ and $Y_{c_2} = 0.$

We claim that $X$ is a submodule of $M$.
For the proof, we use that the maps $\alpha\:X_{c_2} \to X_{c_1},$
$\alpha\:X_{c_1} \to X_{a_0},$
$\alpha\:X_{a_1} \to X_{d_0},$ and
$\beta\:X_{c_2} \to X_{b_1}$  are surjective. Since
$M$ has no direct summand of the form $P(c_2)$, we have $\alpha^3M_{c_2} =
\beta^2M_{c_2} = 0,$ thus $\alpha X_{a_0} = 0$ and $\beta X_{b_1} = 0$.
The relations $\alpha\beta = 0 = \beta\alpha$ show that also the subspaces
$\beta X_{d_0},\ \beta X_{a_0},\ \beta X_{c_1},\ \alpha X_{b_1}$
all are zero.

Let $M'$ be defined by $M'|U = Y$ and $M'_x = M_x$ for those vertices $x$
in the quiver of $\Lambda'_1$ which do not belong to $U$. Clearly, $M'$
is a submodule of $M$ and we have $M = X\oplus M'.$ By construction, $M'_{c_2} = 0$,
thus $M'$ is a $\Lambda_1$-module.
$\s$
	\medskip
{\bf Corollary.} {\it If $M$ is a $\Lambda_2$-module of finite projective dimension,
then $\Omega M$ is a $\Lambda_1$-module.}
	\medskip
Proof. The syzygy-module
$\Omega M$ is a $\Lambda_1'$-module of finite projective dimension, thus
according to Lemma 1 and Lemma 2 a $\Lambda_1$-module. $\s$
	\medskip
{\bf Proposition.} {\it Any $\Lambda_m$-module of finite projective
dimension has projective dimension at most $r+m.$}
	\medskip
Proof. Let $M$ be a $\Lambda_m$-module of finite projective dimension.

First, let $m = 0.$ The algebra $\Lambda_0$ is the product of an algebra
of global dimension $r$ (with vertices $d_0,\dots,d_r$) and an
algebra (with vertices $a_0,c_0,b_0,u,v,w,$ $b_{-1},c_{-1}$) whose
non-projective modules have infinite projective dimension. Thus $\pd M \le r.$

Now, let $m \ge 1.$ Then $\Omega M$ is a $\Lambda_{m-1}$-module of finite
projective dimension. By induction $\pd \Omega M \le r+m-1$, thus $\pd M \le r+m.$
$\s$
	\medskip
{\bf 4. Direct limits.}
The abstract of [HZ] claims that there exist infinite dimensional
$\Lambda$-modules of finite projective dimension which are not
direct limits of finitely generated representations of finite
projective dimension. Apparently, the author refers to the $\Lambda_m$-modules
labelled $M_m$ which are presented in Claim 4 of [HZ]
(these are the only infinite dimensional $\Lambda$-modules
exhibited in the paper; they are used  in order to show that
$\Findim \Lambda_m \ge r+m$).
Indeed, these modules $M_m$
have finite projective dimension, namely $\pd M_m = r+m$.
However, {\it the modules $M_m$ {\bf are} direct limits of finitely
generated modules of finite projective dimension,} as we will show.

For $m\ge 0$ and $t\ge 1$, let us introduce a $\Lambda_m$-module $Z_m[t]$
such that $Z_m[1] = Z_m$, with a submodule $U_{mt} \subset Z_m[t]$, as well as a
map $\phi_{mt}\:Z_m[t] \to Z_{m}[t\!+\!1]$
with kernel $U_{mt}$. Below, we display the modules $Z_m[t]$ with $t = 3$.
The submodule $U_{mt}$ is the zero module in case $m \ge 3$,
and is the shaded part in case $m \le 2$.
The module $X = Z_m[t]/U_{mt}$ has a filtration
$0 \subseteq X_0 \subset X_1 \subset \cdots \subset X_t \subseteq Z_m[t]/U_{mt}$
with isomorphic subfactors
$X_{s}/X_{s-1}$ for $1\le s \le t$.
In our display, we
enclose the subfactors $X_{s}/X_{s-1}$ with $1\le s \le 3$ by dotted lines.
$$
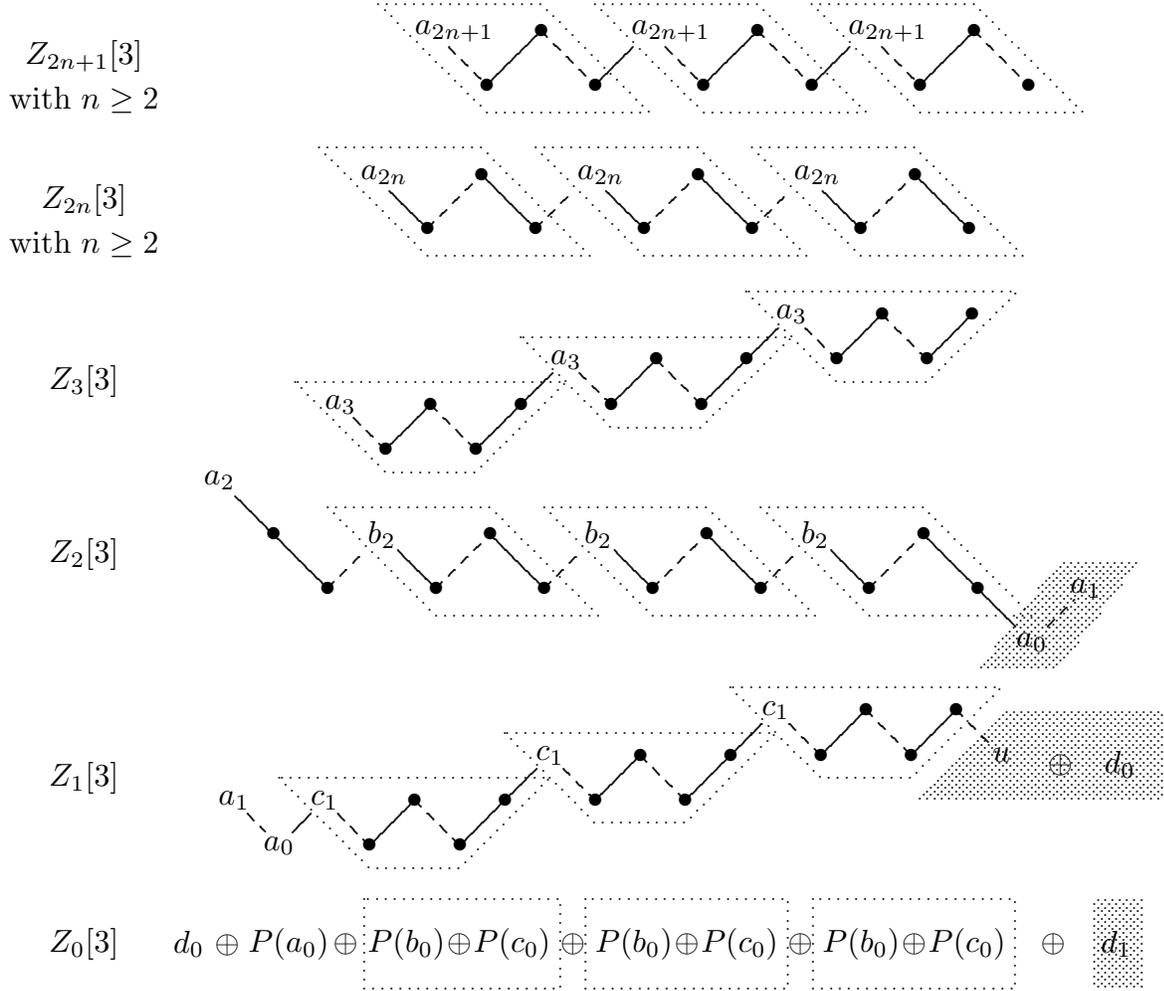
\hbox{\beginpicture
    \setcoordinatesystem units <.9cm,.8cm>
\put{$Z_{2n+1}[3]$} at -6 -0
\put{with $n\ge 2$} at -6 -0.6
\put{$Z_{2n}[3]$} at -6 -2
\put{with $n\ge 2$} at -6 -2.6
\put{$Z_3[3]$} at -6 -4.5
\put{$Z_2[3]$} at -6 -6.9
\put{$Z_1[3]$} at -6 -10
\put{$Z_0[3]$} at -6 -12.3
\put{\beginpicture
    \setcoordinatesystem units <.6cm,.6cm>
\multiput{$\bullet$} at 1 0  2 1  3 0
   5 0  6 1  7 0  9 0  10 1  11 0  /
\multiput{$a_{2n+1}$} at  0.4 1  4.4 1  8.4 1 /
\plot 1 0  2 1 /
\plot 3 0  3.7 0.7 /
\plot 5 0  6 1 /
\plot 7 0  7.7  0.7 /
\plot 9 0  10 1 /

\setdashes <1mm>
\plot 0.3 0.7  1 0 /
\plot 2 1  3 0 /

\plot 4.3 0.7  5 0 /
\plot 6 1  7 0 /

\plot 8.3 0.7  9 0 /
\plot 10 1  11 0 /

\setdots <1mm>
\plot -1 1.5  1 -.5  4 -.5  2 1.5  -1 1.5 /
\plot 3 1.5  5 -.5  8 -.5  6 1.5  3 1.5 /
\plot 7 1.5  9 -.5  12 -.5  10 1.5  7 1.5 /
\endpicture} at 1.9 -0

\put{\beginpicture
    \setcoordinatesystem units <.6cm,.6cm>
\multiput{$\bullet$} at 1 0  2 1  3 0
   5 0  6 1  7 0  9 0  10 1  11 0  /
\multiput{$a_{2n}$} at  0.2 1  4.2 1  8.2 1 /
\setdashes <1mm>
\plot 1 0  2 1 /
\plot 3 0  3.7 0.7 /
\plot 5 0  6 1 /
\plot 7 0  7.7  0.7 /
\plot 9 0  10 1 /

\setsolid
\plot 0.3 0.7  1 0 /
\plot 2 1  3 0 /

\plot 4.3 0.7  5 0 /
\plot 6 1  7 0 /

\plot 8.3 0.7  9 0 /
\plot 10 1  11 0 /

\setdots <1mm>
\plot -1 1.5  1 -.5  4 -.5  2 1.5  -1 1.5 /
\plot 3 1.5  5 -.5  8 -.5  6 1.5  3 1.5 /
\plot 7 1.5  9 -.5  12 -.5  10 1.5  7 1.5 /
\endpicture} at 1.2 -2

\put{\beginpicture
    \setcoordinatesystem units <.5cm,.5cm>
\multiput{$a_3$} at  0 1  5 2  10 3 /
\multiput{$\bullet$} at 1 0  2 1  3 0  4 1
                        6 1  7 2  8 1  9 2
                        11 2  12 3  13 2  14 3 /
\plot 1 0  2 1 /
\plot 3 0  4 1 /
\plot 6 1  7 2 /
\plot 8 1  9 2 /
\plot 11 2  12 3 /
\plot 13 2  14 3 /
\plot 4 1  4.7 1.7 /
\plot 9 2  9.7 2.7 /
\setdashes <1mm>
\plot 0.3  0.7  1 0 /
\plot 2 1  3 0 /
\plot 5.3  1.7  6 1 /
\plot 7 2  8 1 /
\plot 10.3  2.7  11 2 /
\plot 12 3  13 2 /

\setdots <1mm>
\plot -1 1.5  1 -.5  3 -.5  5 1.5  -1 1.5 /
\plot 4 2.5  6 .5  8 .5  10 2.5  4 2.5 /
\plot 9 3.5  11 1.5  13 1.5  15 3.5  9 3.5 /

\endpicture} at 1 -4.5

\put{\beginpicture
    \setcoordinatesystem units <.6cm,.6cm>
\multiput{$\bullet$} at -2 1  -1 0  1 0  2 1  3 0
   5 0  6 1  7 0  9 0  10 1  11 0   /
\put{$a_2$} at -3 2
\put{$a_0$} at 12 -1
\put{$a_1$} at 13 0
\multiput{$b_2$} at  0 1  4 1  8 1 /
\setdashes <1mm>
\plot -1 0  -.3 .7 /
\plot 1 0  2 1 /
\plot 3 0  3.7 0.7 /
\plot 5 0  6 1 /
\plot 7 0  7.7  0.7 /
\plot 9 0  10 1 /
\plot 12.3 -0.7  12.8 -.2 /

\setsolid
\plot -2.7 1.7  -1 0 /
\plot 0.3 0.7  1 0 /
\plot 2 1  3 0 /

\plot 4.3 0.7  5 0 /
\plot 6 1  7 0 /

\plot 8.3 0.7  9 0 /
\plot 10 1  11.7 -.7 /

\setdots <1mm>
\plot -1 1.5  1 -.5  4 -.5  2 1.5  -1 1.5 /
\plot 3 1.5  5 -.5  8 -.5  6 1.5  3 1.5 /
\plot 7 1.5  9 -.5  12 -.5  10 1.5  7 1.5 /

\setshadegrid span <.4mm>
\vshade 11 -1.5 -1.5 <z,z,,> 12.5 -1.5 0.5 <z,z,,> 14 .5 .5  /

\endpicture} at 1 -7

\put{\beginpicture
    \setcoordinatesystem units <.5cm,.5cm>
\multiput{$c_1$} at  0 1  5 2  10 3 /
\put{$u$} at 15 2
\put{$a_1$} at -2 1
\put{$a_0$} at -1 0
\multiput{$\bullet$} at 1 0  2 1  3 0  4 1
                        6 1  7 2  8 1  9 2
                        11 2  12 3  13 2  14 3 /
\plot -.7 0.3  -.3 0.7 /
\plot 1 0  2 1 /
\plot 3 0  4 1 /
\plot 6 1  7 2 /
\plot 8 1  9 2 /
\plot 11 2  12 3 /
\plot 13 2  14 3 /
\plot 4 1  4.7 1.7 /
\plot 9 2  9.7 2.7 /
\setdashes <1mm>
\plot -1.7 0.7  -1.3 0.3 /
\plot 0.3  0.7  1 0 /
\plot 2 1  3 0 /
\plot 5.3  1.7  6 1 /
\plot 7 2  8 1 /
\plot 10.3  2.7  11 2 /
\plot 12 3  13 2 /
\plot 14 3  14.8 2.2 /
\setdots <1mm>
\plot -1 1.5  1 -.5  3 -.5  5 1.5  -1 1.5 /
\plot 4 2.5  6 .5  8 .5  10 2.5  4 2.5 /
\plot 9 3.5  11 1.5  13 1.5  15 3.5  9 3.5 /
\put{$\oplus \quad d_0$} at 17 1.8

\setshadegrid span <.4mm>
\vshade 13 1 1 <z,z,,> 15 1 3 <z,z,,> 18.7 1 3  /

\endpicture} at 1.3 -10

\put{\beginpicture
    \setcoordinatesystem units <.7cm,1cm>
\put{$d_0\ \oplus$} at -2 0
\put{$P(a_0)$} at -.7 0
\put{$d_1$} at 12.4 0
\multiput{$\oplus$} at 0.2 0  3.8 0   7.4 0  11.4 0  /

\multiput{$P(b_0)\!\oplus\! P(c_0)$} at 2 0  5.6 0  9.2 0 /

\setdots <1mm>
\plot .5 .5  .5 -.5  3.6 -.5  3.6 .5  .5 .5 /
\plot 4 .5  4 -.5  7.2 -.5  7.2 .5  4 .5 /
\plot 7.6 .5  7.6 -.5  10.8 -.5  10.8 .5  7.6 .5 /

\setshadegrid span <.4mm>
\vshade 12 -.5 .5 <z,z,,> 12.8 -.5 .5  /

\endpicture} at 1 -12.3

\endpicture}
$$
	\bigskip
The map $\phi_{mt}$ is given by the
obvious embedding of $Z_m[t]/U_{mt}$ into $Z_m[t\!+\!1]$ and
we define $M_m = \lim_{t} (Z_m[t],\phi_{mt}).$
For $m\ge 1$, the modules $M_m$ are those presented in [HZ], Claim 4.
As in section 2, one easily checks that $\Omega (Z_{m+1}[t]) = Z_m[t]$ for
any $m\ge 0$ and $t\ge 1$, so that $\pd Z_m[t] = r+m.$ Also, $\Omega M_{m+1} = M_m,$
and therefore $\pd M_m = r+m.$

{\bf Remark.} For $m\ge 3$, the module $M_m$ is just a Pr\"ufer module
for its support algebra (which is hereditary).
	\bigskip
{\bf Acknowledgment.} The author declares that there is no conflict of interest.
All data generated or analysed during this study are included in this article.

	\bigskip
{\bf 5. Reference.}
\medskip
\item{[HZ]} B\. Huisgen-Zimmermann: Representation-tame algebras need not be homologically
tame. Algebras and Representation Theory 19 (2016), 943--956.
	\medskip
{\baselineskip=1pt
\rmk
C. M. Ringel\par
Fakult\"at f\"ur Mathematik, Universit\"at Bielefeld \par
POBox 100131, D-33501 Bielefeld, Germany  \par
ringel\@math.uni-bielefeld.de\par}
\smallskip

\vfill\eject
{\bf Appendix. The shape of the indecomposable projective $\Lambda_5$-modules.}
These graphical displays can be found in [HZ]. But the referee has suggested to
provide the pictures also here.
$$
\hbox{\beginpicture
    \setcoordinatesystem units <.65cm,.65cm>
\put{\beginpicture
\put{$d_{2i}$} at 0 0
\put{} at -1 -1
\put{$d_{2i+1}$} at 1 -1
\setdashes <.8mm>
\plot .3 -.3  .75 -.75 /
\put{for $0 \le 2i < r$} at 3.3 -.5
\endpicture} at 0 8
\put{\beginpicture
\put{$d_{2i+1}$} at 0 0
\put{$d_{2i+2}$} at -1 -1
\put{} at 1 -1
\plot -.3 -.3  -.8 -.8 /
\setdashes <.8mm>
\put{for $0 < 2i+1 < r$} at 3 -.5
\endpicture} at 7 8
\put{$d_r$} at 11.5 8.55
\put{\beginpicture
\put{$u$} at 0 0
\put{$u$} at -1 -1
\put{} at 1 -1
\plot -.2 -.2  -.8 -.8 /
\setdashes <.8mm>
\endpicture} at 0 6
\put{\beginpicture
\put{$v$} at 0 -1
\put{} at -1 -1
\put{$v$} at 1 0
\plot .2 -.8  .8 -.2 /
\endpicture} at 2 6
\put{\beginpicture
\put{$w$} at 0 0
\put{$w$} at -1 -1
\put{} at 1 -1
\plot -.2 -.2  -.8 -.8 /
\setdashes <.8mm>
\endpicture} at 5.3 6
\put{\beginpicture
\put{$b_{-1}$} at 0 0
\put{$b_{-1}$} at 1 -1
\put{} at -1 -1
\setdashes <.8mm>
\plot .35 -.35  .75 -.75 /
\endpicture} at 8 6

\put{\beginpicture
\put{$c_{-1}$} at 0 0
\put{} at -1 -1
\put{$c_{-1}$} at 1 -1
\setdashes <.8mm>
\plot .35 -.35  .75 -.75 /
\endpicture} at 10.5 6
\put{\beginpicture
\put{$a_0$} at 0 0
\put{$c_0$} at -1 -1
\put{$u$} at 1 -1
\put{$c_{-1}$} at -1 -2
\put{\strut} at -1 -3
\plot -.2 -.2  -.8 -.8 /
\plot -1 -1.3  -1 -1.7 /
\setdashes <.8mm>
\plot .3 -.3  .75 -.75 /
\endpicture} at 0 3
\put{\beginpicture
\put{$b_0$} at 0 0
\put{$b_{-1}$} at -1 -1
\put{$v$} at 1 -1
\put{\strut} at -1 -2
\put{\strut} at -1 -3
\plot -.2 -.2  -.8 -.8 /
\setdashes <.8mm>
\plot .3 -.3  .75 -.75 /
\endpicture} at 4 3
\put{\beginpicture
\put{$c_0$} at 0 0
\put{$c_{-1}$} at -1 -1
\put{$w$} at 1 -1
\put{\strut} at -1 -3
\plot -.2 -.2  -.8 -.8 /
\setdashes <.8mm>
\plot .3 -.3  .75 -.75 /
\endpicture} at 8 3

\put{\beginpicture
\put{$a_1$} at 0 0
\put{$d_0$} at -1 -1
\put{$a_0$} at 1 -1
\put{$u$} at 1 -2
\put{\strut} at -1 -3
\plot -.2 -.2  -.8 -.8 /
\setdashes <.8mm>
\plot .3 -.3  .75 -.75 /
\plot 1 -1.3  1 -1.9 /
\endpicture} at 0 0
\put{\beginpicture
\put{$b_1$} at 0 0
\put{$b_0$} at -1 -1
\put{$c_0$} at 1 -1
\put{$b_{-1}$} at -1 -2
\put{$w$} at 1 -2
\put{\strut} at -1 -3
\plot -.2 -.2  -.8 -.8 /
\plot -1 -1.3  -1 -1.7 /
\setdashes <.8mm>
\plot .3 -.3  .75 -.75 /
\plot 1 -1.3  1 -1.9 /
\endpicture} at 4 0
\put{\beginpicture
\put{$c_1$} at 0 0
\put{$a_0$} at -1 -1
\put{$b_0$} at 1 -1
\put{$c_0$} at -1 -2
\put{$v$} at 1 -2
\put{$c_{-1}$} at -1 -3
\plot -.2 -.2  -.8 -.8 /
\plot -1 -1.3  -1 -1.7 /
\plot -1 -2.3  -1 -2.7 /
\setdashes <.8mm>
\plot .3 -.3  .75 -.75 /
\plot 1 -1.3  1 -1.9 /
\endpicture} at 8 0

\put{\beginpicture
\put{$a_2$} at 0 0
\put{$c_2$} at -1 -1
\put{$a_1$} at 1 -1.5
\put{$c_1$} at -1 -2
\put{$a_0$} at 0 -3
\plot -.2 -.2  -.8 -.8 /
\plot -1 -1.3  -1 -1.7 /
\plot -.8 -2.2  -.2 -2.8 /
\setdashes <.8mm>
\plot .2 -.3  .8 -1.3 /
\plot .2 -2.8  .8 -1.7 /
\endpicture} at 0 -4
\put{\beginpicture
\put{$b_2$} at 0 0
\put{$b_1$} at -1 -1
\put{$c_1$} at 1 -1
\put{$b_0$} at 0 -2
\put{\strut} at 0 -3
\plot -.2 -.2  -.8 -.8 /
\plot -.8 -1.2  -.2 -1.8 /
\setdashes <.8mm>
\plot .2 -.3  .75 -.75 /
\plot .75 -1.25  .2 -1.8 /
\endpicture} at 4 -4
\put{\beginpicture
\put{$c_2$} at 0 0
\put{$c_1$} at -1 -1
\put{$b_1$} at 1 -1.5
\put{$a_0$} at -1 -2
\put{$c_0$} at 0 -3
\plot -.2 -.2  -.8 -.8 /
\plot -1 -1.3  -1 -1.7 /
\plot -.8 -2.2  -.2 -2.8 /
\setdashes <.8mm>
\plot .2 -.3  .8 -1.3 /
\plot .2 -2.8  .8 -1.7 /
\endpicture} at 8 -4
\put{\beginpicture
\put{$a_3$} at 0 0
\put{$a_2$} at -1 -1
\put{$b_2$} at 1 -1.5
\put{$c_2$} at -1 -2
\put{$c_1$} at 0 -3
\plot -.2 -.2  -.8 -.8 /
\plot -1 -1.3  -1 -1.7 /
\plot -.8 -2.2  -.2 -2.8 /
\setdashes <.8mm>
\plot .2 -.3  .8 -1.3 /
\plot .2 -2.8  .8 -1.7 /
\endpicture} at 0 -8

\put{\beginpicture
\put{$b_3$} at 0 0
\put{$b_2$} at -1 -1
\put{$c_2$} at 1 -1
\put{$b_1$} at 0 -2
\put{\strut} at 0 -3
\plot -.2 -.2  -.8 -.8 /
\plot -.8 -1.2  -.2 -1.8 /
\setdashes <.8mm>
\plot .2 -.3  .75 -.75 /
\plot .75 -1.25  .2 -1.8 /
\endpicture} at 4 -8

\put{\beginpicture
\put{$a_4$} at 0 0
\put{$b_3$} at -1 -1
\put{$a_3$} at 1 -1
\put{$b_2$} at 0 -2
\put{\strut} at 0 -3
\plot -.2 -.2  -.8 -.8 /
\plot -.8 -1.2  -.2 -1.8 /
\setdashes <.8mm>
\plot .2 -.3  .75 -.75 /
\plot .75 -1.25  .2 -1.8 /
\endpicture} at 0 -12
\put{\beginpicture
\put{$b_4$} at 0 0
\put{$a_3$} at -1 -1
\put{$b_3$} at 1 -1.5
\put{$a_2$} at -1 -2
\put{$c_2$} at 0 -3
\plot -.2 -.2  -.8 -.8 /
\plot -1 -1.3  -1 -1.7 /
\plot -.8 -2.2  -.2 -2.8 /
\setdashes <.8mm>
\plot .2 -.3  .8 -1.3 /
\plot .2 -2.8  .8 -1.7 /
\endpicture} at 4 -12

\put{\beginpicture
\put{$a_5$} at 0 0
\put{$a_4$} at -1 -1
\put{$b_4$} at 1 -1
\put{$b_3$} at 0 -2
\put{} at 0 -3
\plot -.2 -.2  -.8 -.8 /
\plot -.8 -1.2  -.2 -1.8 /
\setdashes <.8mm>
\plot .2 -.3  .75 -.75 /
\plot .75 -1.25  .2 -1.8 /
\endpicture} at 0 -16
\put{\beginpicture
\put{$b_5$} at 0 0
\put{$b_4$} at -1 -1
\put{$a_4$} at 1 -1
\put{$a_3$} at 0 -2
\put{} at 0 -3
\plot -.2 -.2  -.8 -.8 /
\plot -.8 -1.2  -.2 -1.8 /
\setdashes <.8mm>
\plot .2 -.3  .75 -.75 /
\plot .75 -1.25  .2 -1.8 /
\endpicture} at 4 -16

\endpicture}
$$

\bye